\documentclass[12pt]{article}
\usepackage[russian,english]{babel}
\usepackage{amsmath}
\usepackage[standard,thmmarks,amsmath]{ntheorem}
\usepackage{amssymb}
\usepackage{graphicx}

\def\le{\leqslant}

\def\ge{\geqslant}

\title{Minimizing Degree-based Topological Indices for Trees with Given Number of Pendent Vertices\thanks{This research is supported by the grant of Russian Foundation for Basic Research, project No 13-07-00389.}
}
\author{Mikhail Goubko\thanks{\texttt{mgoubko@mail.ru}}\\\emph{Institute of Control Sciences of RAS}}

\begin{document}

\maketitle

\begin{abstract}
We derive sharp lower bounds for the first and the second Zagreb
indices ($M_1$ and $M_2$ respectively) for trees with the given
number of pendent vertices and find optimal trees. $M_1$ is
minimized by a tree with all internal vertices having degree $4$,
while $M_2$ is minimized by a tree where each ``stem'' vertex is
incident to $3$ or $4$ pendent vertices and one internal vertex,
while the rest internal vertices are incident to 3 other internal
vertices. The technique is shown to generalize to the weighted first
Zagreb index, the zeroth order general Randi\'{c} index, as long as
to many other degree-based indices.
\end{abstract}

\section*{Introduction}

Topological graph indices are widely used in mathematical chemistry
to predict properties of chemical compounds. They have been
intensively studied in recent years. Dozens of various indices were
suggested \cite{Balaban} to describe topology of complex molecules,
among the earliest and the most famous being the first and the
second Zagreb indices -- $M_1$ and $M_2$ respectively
\cite{GutmanTrin72}. The popular research problem is to find lower
and upper bounds of an index over a certain set of graphs and to
characterize extremal graphs in this set.

The typical set to study is that of all graphs (trees, bipartite or
unicyclic graphs, ``cacti'', etc) of the fixed order (i.e. with the
fixed number of vertices). Extremal graphs on these sets often
appear to be degenerate. For example, the chain minimizes Zagreb
indices, while the star maximizes them over the set of trees of
order $N$ (see \cite{DasGutman04,GutmanDas04}). Even when the set of
admissible graphs is cut (by limiting degrees, chromatic or matching
numbers, etc), extremal graphs are typically found on the
``boundary'' of the set. For instance, the ``broom'' (i.e., the star
$K_{1,\Delta}$ with the path of length $N-\Delta-1$ attached to any
pendent vertex) minimizes $M_2$ over the set of trees with the fixed
maximum degree $\Delta$ (see \cite{Stevanovich12}), the path of
length $N-k$ attached to the cycle of length $k$ minimizes both
$M_1$ and $M_2$ over the set of all unicyclic graphs of order $N$
and girth $k$ (see \cite{Deng07}), etc.

We optimize indices over the set of trees with the fixed number of
pendent vertices. If hydrogen atoms are not suppressed from Sachs
diagrams \cite{Polansky75}, this set can be interpreted as that of
all acyclic molecules with the fixed number of hydrogen atoms. In
hydrogen-suppressed diagrams of paraffins pendent vertices stand for
methyl groups $\text{CH}_3$.

This set of graphs is of interest as it provides a ``vertex-number
vs degree'' trade-off for degree-based indices, resulting in
optimality of nontrivial internal solutions. Note that such
``internal'' solutions do not arise even when one studies the set of
graphs parameterized by the number of pendent vertices $n$ and the
total number of vertices $N$. For example, the star $K_{1,n}$ with
$n$ (roughly equal) paths attached to its rays maximizes $M_2$ over
the set of trees with fixed $n$ and $N$ \cite{LiuLiu12}. The root of
the star in this graph has the maximum possible degree $n$ while all
other internal vertices have the minimum possible degree 2.
Unicyclic graphs with minimum possible vertex degrees (no more than
3) minimize both $M_1$ and $M_2$ over the set of ``cacti'' with
fixed $n$ and $N$ \cite{Lee12}. For more results on extremal trees
with fixed $N$ and $n$ for the Randi\'{c} index (which is closely
related to $M_2$) one can refer to \cite{LiShi08}.

Below we show that in the tree minimizing $M_1$ over the set of all
trees with $n$ pendent vertices almost all internal vertices have
degree $4$, which is strictly greater than the minimum possible
degree $3$ but less than the maximum possible degree $n$. We also
show that in a tree, which minimizes $M_2$, internal vertices have
degrees $3$, $4$ and $5$. Even more surprising structures are shown
to minimize the generalized Randi\'{c} index or the multiplicative
Zagreb indices $\Pi_1$ and $\Pi_2$.

\section{The First Zagreb Index}

Let $G$ be a simple connected undirected graph with the vertex set
$V(G)$ and the edge set $E(G)$. Denote by $d_G(v)$ the degree of a
vertex $v\in V(G)$ in the graph $G$, i.e., the number of vertices
being incident to $v$ in $G$. The first Zagreb index is defined in
\cite{GutmanTrin72} as
\begin {equation}\label{M1_def}
M_1(G):=\sum_{v\in V(G)} d_G(v)^2,
\end{equation}
while the second Zagreb index -- as
\begin {equation}\label{M2_def}
M_2(G):=\sum_{uv\in E(G)} d_G(u)d_G(v).
\end{equation}

The vertex $v\in V(G)$ with $d_G(v)=1$ is called a pendent vertex.
Denote the set of pendent vertices of the graph $G$ with $W(G)$. A
connected graph $T$ with $N$ vertices and $N-1$ edges is called a
tree.

\begin{theorem}\label{LB_M1_theorem} For any tree $T$ with $n\ge 2$ pendent vertices
$M_1(T)\ge 9n-16$ if $n$ is even. The equality holds if $T$ is a
$4$-tree \emph{(}with $d_T(m)=4$ for all $m \in V(G)\backslash
W(G)$\emph{)}. If $n$ is odd, then $M_1(T)\ge 9n-15$, and the
equality holds if $T$ is a tree with all internal vertices having
degree $4$ except the one of degree $3$.
\end{theorem}
\begin{proof}
For $n=2$ the optimal tree is the complete graph $K_2$, and the
theorem obviously holds. If $n>2$, there must be at least one
internal vertex in a tree.

Note that the tree $T$ cannot minimize $M_1(T)$ over the set of all
trees with $n$ vertices if it contains an internal vertex of degree
$2$. Actually, the index is reduced by eliminating such a vertex and
shortcuting its incident vertices. So, below we restrict attention
to the trees with internal vertex degrees at least $3$.

For an arbitrary tree with $n>2$ pendent vertices and $q>0$ internal
vertices of degrees $d_1, ..., d_q$ the following identity holds:
\begin{equation}\label{degrees_identity}
n-2=\sum_{i=1}^q (d_i-2).
\end{equation}
Thus, minimization of $M_1$ for fixed $n$ and $q$ reduces to
minimization of $n 1^2+\sum_{i=1}^q d_i^2$ over all $d_i = 3,4,
...$, $i=1,...,q$ satisfying (\ref{degrees_identity}). Ignoring
integer constraints from the first order conditions we obtain an
obvious solution of this convex program: $d_i=2+(n-2)/q$ for all
$i=1,...,q$. Then, to find optimal $q$ we minimize
$n+q(2+(n-2)/q)^2$ over all $q=1,...,n-2$ (the range follows from
(\ref{degrees_identity})). Relaxing the integer constraint from the
first order condition find optimal $q=(n-2)/2$ and $d_i=4$. Thus, as
we relaxed some integer constraints during minimization, $M_1(T)\ge
n+16(n-2)/2=9n-16$. It follows from (\ref{degrees_identity}) that
the tree $T_{even}$ with $q=(n-2)/2$ internal vertices of degree $4$
exists for even $n$. An obvious calculation gives
$M_1(T_{even})=9n-16$. For odd $n$ it follows from
(\ref{degrees_identity}) that no $4$-tree exists and, thus, the
lower-bound estimate $9n-16$ cannot be achieved. At the same time,
there exists a tree $T_{odd}$ with all internal vertices having
degree $4$ except the one of degree $3$ with $M_1(T_{odd})=9n-15$.
As the index $M_1$ is integer-valued, $T_{odd}$ is optimal for odd
$n$.
\end{proof}

The above theorem says that, at least in the considered stylized
setting, carbon of valency $4$ is the best connector for any given
number of hydrogen atoms (if hydrogen atoms are not suppressed from
the molecular graph) or methyl groups $\text{CH}_3$ (in
hydrogen-suppressed diagrams) in terms of minimization of the first
Zagreb index. In both cases $M_1$ is minimized by alkanes
$\text{C}_m\text{H}_{2m+2}$.

Let us account for heterogeneity of atoms by adding to every term
$d_G(v)^2$ in $M_1(G)$ a weight depending on the vertex degree (the
valency of an atom in a molecule). The following theorem gives the
lower bound for the generalized index $C(G):=\sum_{v\in
V(G)}c(d_G(v))$, where $c(\cdot)$ is an arbitrary non-negative
function of the vertex degree. As we minimize $C(G)$, it is natural
to call it the cost of the graph $G$, and to call $c(d_G(m))$ the
cost of the vertex $m$ in the graph $G$.

\begin{theorem}\label{LB_M1_general_theorem}
For $n \ge 2$
\begin {equation}\label{LB_M1_general}
C(T)\ge
\underline{C}(n):=nc(1)+(n-2)\frac{c(\Delta(n))}{\Delta(n)-2}
\end{equation}
for an arbitrary tree $T$ with $n$ pendent vertices, where
\begin{equation}\label{optimal_degree}
\Delta(2)=3, \Delta(n) \in \text{\normalfont Argmin}_{d=3,...,n}
\frac{c(d)}{d-2}\text{ \normalfont for }n>2.
\end{equation}
When $q(n):=\frac{n-2}{\Delta(n)-2}$ is integer, the equality in
\emph{(\ref{LB_M1_general})} is achieved at an arbitrary tree, where
all internal vertices have degree $\Delta(n)$.
\end{theorem}
\begin{proof}
Fix an arbitrary pendent vertex $w\in W(T)$ in a tree $T$. Then
$C(T)=c(1)+\sum_{v\in V(T) \backslash\{w\}}c(d_T(v))$. Let us call
the tree $T$ with the selected pendent vertex $w$ \emph{the attached
tree with the root $w$} and define the cost of this attached tree as
$C_a(T,w):=C(T)-c(1)$. So, the root is still a pendent vertex of an
attached tree, but the cost of the root is not included in the cost
of the attached tree. We will also refer to the vertex incident to
the root in an attached tree as to the ``\emph{sub-root}''.

The set of trees with $n$ pendent vertices coincides with that of
attached trees with $n$ pendent vertices, and their costs differ
only by a constant. So, the problem of cost minimization for a tree
with $n$ pendent vertices is equivalent to cost minimization for an
attached tree with $n$ pendent vertices. Below we prove by induction
that for any attached tree $T$ with $n \ge 2$ pendent vertices
\begin{equation}\label{attached_inequality}
C_a(T, \cdot) \ge
\underline{C}_a(n):=(n-1)c(1)+(n-2)\frac{c(\Delta(n))}{\Delta(n)-2}.
\end{equation}

For $n=2$ (\ref{attached_inequality}) is satisfied as equality, as
the optimal attached tree is a complete graph $K_2$ with $C_a(K_2,
\cdot)=c(1)$ (remember the cost of the root is not counted). Suppose
(\ref{attached_inequality}) is valid for all $n'<n$. Let us prove
that it is also valid for any attached tree $T$ with $n$ pendent
vertices and some root $w$.

As $n \ge 3$, the sub-root $m$ of $T$ is an internal vertex. So, the
cost $C_a(T,w)$ can be written as the sum of the cost of the
sub-root $m$ and costs of the sub-trees $T_1, ..., T_{d_T(m)-1}$
with $n_1, ..., n_{d_T(m)-1}$ pendent vertices respectively,
attached to $m$:
$$C_a(T,w)=c(d_T(m))+\sum_{i=1}^{d_T(m)-1}C_a(T_i, m).$$

As $n_i<n$, by induction hypothesis
$$C_a(T,w) \ge c(d_T(m))+\sum_{i=1}^{d_T(m)-1}\left((n_i-1)c(1)+(n_i-2)\frac{c(\Delta(n_i))}{\Delta(n_i)-2}\right).$$
Note that from (\ref{optimal_degree}) follows that
$$\frac{c(\Delta(n_i))}{\Delta(n_i)-2}\ge \frac{c(\Delta(n))}{\Delta(n)-2},$$
and also that $n_1+...+n_{d_T(m)-1}=n+d_T(m)-2$. Thus,
\begin{eqnarray}
C_a(T,w) \ge c(d_T(m))+\sum_{i=1}^{d_T(m)-1}\left\{(n_i-1)
c(1)+(n_i-2)\frac{c(\Delta(n))}{\Delta(n)-2}\right\}=\nonumber\\
=(n-1)c(1)+c(d_T(m))+(n-d_T(m))\frac{c(\Delta(n))}{\Delta(n)-2}.
\end{eqnarray}
Obviously, $3 \le d_T(m)\le n$, so
\begin{eqnarray}
C_a(T,w) \ge
(n-1)c(1)+\min_{d=3,...,n}\left\{c(d)+(n-d)\frac{c(\Delta(n))}{\Delta(n)-2}\right\}=\nonumber\\
(n-1)c(1)+(n-2)\frac{c(\Delta(n))}{\Delta(n)-2}+\min_{d=3,...,n}(d-2)\left[\frac{c(d)}{d-2}-\frac{c(\Delta(n))}{\Delta(n)-2}\right].
\end{eqnarray}

From (\ref{optimal_degree}) we know that the expression in square
brackets achieves its minimum (which is equal to zero) at
$d=\Delta(n)$. So, the minimum of the product
$(d-2)\left[\frac{c(d)}{d-2}-\frac{c(\Delta(n))}{\Delta(n)-2}\right]$
is also zero, and
$$C_a(T,w) \ge (n-1)c(1)+(n-2)\frac{c(\Delta(n)}{\Delta(n)-2}.$$
So, inequality (\ref{LB_M1_general}) is proved.

When $q(n)=\frac{n-2}{\Delta(n)-2}$ is integer, there exists a
$\Delta(n)$-tree $T^*$ with $n$ pendent and $q(n)$ internal
vertices, which has the cost
$C(T^*)=nc(1)+q(n)c(\Delta(n))=nc(1)+(n-2)\frac{c(\Delta(n))}{\Delta(n)-2}=\underline{C}(n).$

This completes the proof\footnote{We use the scheme of the proof
from \cite{Goubko08}, where the similar result was obtained for
directed trees under a more general cost function.}.
\end{proof}

\begin{example}
The above theorem covers the first Zagreb index with $c(k)=k^2$ and
the zeroth order general Randi\'{c} index $c(k)=k^\alpha$ as special
cases. In particular, using (\ref{LB_M1_general}) one can show that
for $\alpha \ge {\ln 2}/{\ln (4/3)}$ a $3$-tree is optimal (and
exists for all $n\ge 2$), while for
$$\alpha \in
\left[{\ln\left(\frac{d-1}{d-2}\right)}/{\ln\left(\frac{d+1}{d}\right)},{\ln\left(\frac{d-2}{d-3}\right)}/{\ln\left(\frac{d}{d-1}\right)}\right)$$
the optimal degree $\Delta(n)=d$, where $d=4,5,...$, for $n\ge d$
(this means that $d$-tree is optimal for $n\ge d$ when such a tree
exists). For $\alpha \le 1$ the optimal tree is a star $K_{1,n}$, as
$\frac{d^\alpha}{d-2}$ in (\ref{optimal_degree}) is monotone
decreasing for $d\ge 3$.
\end{example}

\begin{example}
The first and the second multiplicative Zagreb indices were defined
in \cite{Gutman11} as
$$\Pi_1(G):=\prod_{v\in V(G)}d_G(v)^2, \Pi_2(G):=\prod_{uv\in E(G)}d_G(u)d_G(v).$$
Instead of summation, as in (\ref{M1_def}) and (\ref{M2_def}),
contributions of vertices here (in the case of the first index) or
edges (in the case of the second index) are multiplied.

Minimization of $\Pi_1(G)$ reduces to minimization of
$$C(G):=\ln\Pi_1(G)=2\sum_{v\in V(G)}\ln d_G(v).$$
From Theorem \ref{LB_M1_general_theorem}, as $\frac{\ln d}{d-2}$ is
monotone decreasing for $d\ge 3$, the tree with $n$ pendent vertices
minimizing $\Pi_1$ is a star $K_{1,n}$.

It is shown in \cite{Gutman11} that for an arbitrary tree $T$
$$\Pi_2(T)=\prod_{v\in V(T)}d_T(v)^{d_T(v)}.$$
So, minimization of $\Pi_2(T)$ for all trees with $n$ pendent
vertices is equivalent to minimization of
$C(T):=\ln\Pi_2(T)=\sum_{v\in V(T)}d_T(v)\ln d_T(v)$. Set $c(d)=d
\ln d$ in (\ref{optimal_degree}) and obtain $\Delta(n) = \min[n,5]$.
Then, from (\ref{LB_M1_general}) we see that $\Pi_2(T)=\exp(C(T))\ge
\exp\left(\frac{5 \ln 5}{3}(n-2)\right)$ for $n\ge 5$ with equality
at any $5$-tree when $(n-2) / 3$ is integer.
\end{example}

When $\frac{n-2}{\Delta(n)-2}$ is not integer, there exists no
$\Delta(n)$-tree with $n$ pendent vertices, and the lower bound
(\ref{LB_M1_general}) is not sharp. Nevertheless, for every specific
function $c(d)$ one often can prove the optimal tree to be a some
minimal perturbation of the $\Delta(n)$-tree. Typically the optimal
tree is a \emph{bidegree tree}, where almost all internal vertices
have degree $\Delta(n)$, while several vertices have degree
$\Delta(n)+1$ or $\Delta(n)-1$ (like in Theorem
\ref{LB_M1_theorem}).

\section{The Second Zagreb Index}

An internal vertex in a tree is called a \emph{stem vertex} if it
has incident pendent vertices (see \cite{ChenLiu12}). The edge
connecting a stem with a pendent vertex will be referred to as a
\emph{stem edge}.

\begin{theorem}\label{M_2_theorem}
For any tree $T$ with $n\ge 8$ pendent vertices $M_2(T)\ge 11n-27$.
The equality holds if each stem vertex in $T$ has degree $4$ or $5$
while other internal vertices having degree $3$. At least one such
tree exits for any $n\ge 9$.
\end{theorem}
\begin{proof}
Let us employ again the idea of an attached tree from Theorem
\ref{LB_M1_general_theorem}. Below we suggest a suitable
generalization of the concept of an attached tree, then we
interrelate its cost with $M_2$, and, finally, use induction on $n$
to prove the lower bound. The cost of trees, which minimize $M_2$,
is found by a direct calculation.

Let us allow the root of an attached tree to have arbitrary degree
$p\ge 1$. Actually we do not add vertices incident to the root -- it
is still incident only to the sub-root -- but the degree of the root
is substituted to the contribution of the edge $wm$ connecting the
root $w$ with the sub-root $m$ to the index $M_2$. For the attached
tree $T$ with the root $w$ of degree $p$ and the sub-root $m$ of
degree $d$ define its \emph{cost} as
\begin{equation}\label{M2_attached_tree_cost}
C_a(T,w,p):=pd+\sum_{uv\in E(T)\backslash \{wm\}}d_T(u)d_T(v).
\end{equation}

We will consider the root as a pendent vertex only when its degree
$p=1$. Note that it implies the following interrelation between
$M_2(T)$ and the cost of the attached tree $T$ with an arbitrary
root $w\in W(T)$: $M_2(T)=C_a(T,w,1)$. So, the problem of
minimization of $M_2$ over the set of all trees with $n$ pendent
vertices is equivalent to the problem of minimization of the cost of
an attached tree with $n$ vertices and the root of degree $1$.

First we use induction to show that for any attached tree $T$ with
$n$ pendent vertices and some root $w$ of degree $p\ge 3$
\begin{equation}\label{M2_attached_inequality}
C_a(T,w,p)\ge \begin{cases} p, & \mbox{if } n=1,
\\ 11n+3p-18, & \mbox{if } n\ge 2. \end{cases}
\end{equation}

Note, that, as before, we can restrict attention to the trees where
all internal vertices (including the root) have degree at least $3$.
For $n=1$ the inequality (\ref{M2_attached_inequality}) trivially
holds as the only attached tree has only one edge. From
(\ref{M2_attached_tree_cost}), its cost is $p$.

Suppose inequality (\ref{M2_attached_inequality}) holds for all
$n'<n$. Let us prove that it also holds for $n$. As $n\ge2$, the
sub-root of any attached tree is an internal vertex. Consider a tree
$T$ with some root $w$ of degree $p\ge3$ and the sub-root $m$ of
degree $d\ge 3$ having $\delta\ge 0$ incident pendent vertices and
$\Delta \ge 0$ incident internal vertices. Note that
$d=\delta+\Delta+1$ and $3\le d \le n+1$. The cost of the attached
tree $T$ consists of the cost $pd$ of the edge $mw$, the total cost
$\delta d$ of $\delta$ pendent vertices being incident to $m$, and
the sum of costs of $\Delta$ sub-trees $T_1, ..., T_{\Delta}$
attached to $m$: $C_a(T,w,p)=pd+\delta d+\sum_{i=1}^\Delta
C_a(T_i,m,d)$.

To estimate $C_a(T,w,p)$ consider separately the case of $\Delta=0$
and that of $\Delta>0$:
\begin{enumerate}
    \item If $\Delta=0$ then $T=K_{1,n+1}$ with $\delta=n$ and $d=n+1$, so
$C_a(K_{1,n+1},w,p)=(p+n)(n+1)$. Denote $C_1:=(p+n)(n+1)$ for short.
    \item Suppose $\Delta\ge 1$ and let the tree $T_i$ have $n_i\ge 2$
    pendent vertices, $i=1,...,\Delta$. As $2\le n_i<n$, by induction hypothesis $C_a(T_i,m,d)\ge
    11n_i+3d-18$. Taking into account the balance equation $\sum_{i=1}^\Delta
    n_i=n-\delta$, we can estimate the cost of the attached tree
    from below:
    \begin{eqnarray}\label{M_2_attached_inequality2}
    C_a(T,w,p)\ge pd+\delta d+11n-11\delta
    +3(d-\delta-1)(d-6)=\nonumber\\
    11n+pd+\delta(7-2d)+3(d-1)(d-6).
    \end{eqnarray}
    As $n_i\ge 2$, it follows that $3\le d\le n$.
    Also, from $\Delta\ge 1$ and from $d=\delta+\Delta+1$ it follows that $0\le \delta \le
    d-2$. Let us find $d$ and $\delta$ which minimize the right-hand
    side (r.h.s.) in (\ref{M_2_attached_inequality2}).
    Below we consider separately the case of $d=3$ and that of $d\ge
    4$:

\begin{itemize}
    \item If $d=3$ then $7-2d>0$, and r.h.s in
    (\ref{M_2_attached_inequality2}) achieves minimum at $\delta=0$ and equals
    $11n+3p+3(d-1)(d-6)$, which reduces to $C_2:=11n+3p-18$.
    \item If $d\ge 4$ then $7-2d<0$, so r.h.s.
    in (\ref{M_2_attached_inequality2}) achieves its minimum at $\delta=d-2$,
    and equals $11n+pd+d^2-10d+4$. For $p\ge
    3$ and $d\ge 4$ this expression is monotone in $d$ and, thus, r.h.s in (\ref{M_2_attached_inequality2})
    is not less than $11n+4p-20$, which is greater than $C_2$ for $p\ge 3$.
\end{itemize}
So, we conclude that if $\Delta\ge 1$, then $C_a(T,w,p)\ge
C_2=11n+3p-18$.
\end{enumerate}
Let us compare cases 1 and 2 and prove that $C_1$ is never less than
$C_2$ for $n\ge 2$ and $p\ge 3$. Actually, the difference
$C_1-C_2=p(n-2)-10n+n^2+18$ is monotone in $p$, and, thus, achieves
its minimum at $p=3$. Substituting $p=3$ we find that $C_1-C_2\ge
n^2-7n+12$, which is non-negative for all integer $n$.

So, we proved inequality (\ref{M2_attached_inequality}). Let us use
it now to prove that for $p=1$ and $n \ge 9$ $C_a(T,w,1)\ge 11n-27$.

\vspace{12pt}

For $n\ge 3$ in the attached tree $T$ with $n$ pendent vertices
(including the root $w$, as $p=1$) the sub-root $m$ is an internal
vertex. Let the sub-root $m$ have degree $d\ge3$ which adds up from
$\delta \ge 1$ incident pendent vertices (including the root) and
$\Delta\ge0$ internal vertices.

The cost of the attached tree $T$ consists of the total cost $\delta
d$ of $\delta$ pendent vertices incident to $m$ and the sum of costs
of $\Delta$ sub-trees $T_1, ..., T_{\Delta}$ attached to the
sub-root $m$: $C_a(T,w,1)=\delta d+\sum_{i=1}^\Delta C_a(T_i,m,d)$.

\begin{enumerate}
    \item If $\Delta=0$, then $T=K_{1,n}$ with $\delta=d=n$, so
$C_a(K_{1,n},w,1)=n^2$.
    \item Suppose $\Delta\ge 1$ and let the tree $T_i$ have $n_i\ge 2$
    pendent vertices, $i=1,...,\Delta$.
    As $n_i\ge 2$ and $d\ge 3$, from (\ref{M2_attached_inequality}) $C_a(T_i,m,d)\ge
    11n_i+3d-18$. Accounting for the balance equalities $\sum_{i=1}^\Delta
    n_i=n-\delta$ and $\delta+\Delta=d$, we estimate the cost of the
    attached tree as
    \begin{eqnarray}\label{M_2_attached_inequality3}
    C_a(T,w,1)\ge \delta d+11n-11\delta
    +3(d-\delta)(d-6)=\nonumber\\
    =11n+\delta(7-2d)+3d(d-6).
    \end{eqnarray}
    As $n_i\ge 2$, it follows that $3\le d\le n-1$.
    Also recall that $1\le \delta \le
    d-1$. Let us minimize r.h.s. in (\ref{M_2_attached_inequality3}) over all $d=\overline{3,n-1}$ and $\delta=\overline{1,d-1}$. The arguments are similar to that in the case of $p\ge3$:

    \begin{itemize}
    \item If $d=3$, then $7-2d>0$ and r.h.s. in
    (\ref{M_2_attached_inequality3}) achieves its minimum $11n-26$ at $\delta=1$.
    \item If $d\ge 4$, then $7-2d<0$ and r.h.s.
    in (\ref{M_2_attached_inequality3}) achieves minimum $11n+(d-1)(7-2d)+3d(d-6)=11n+d^2-9d-7$ at $\delta=d-1$. Minimum of $11n+d^2-9d-7$ over all integer $d$ is achieved at
    $d=4$ and $d=5$ and is equal to $11n-27$. This is one less than $11n-26$ which we had in the previous case of
    $d=3$.\footnote{This point in the proof is mentioned below in the discussion as a clue to the result for chemical graphs.}
    \end{itemize}
So, finally we conclude that if $\Delta\ge 1$, then $C_a(T,w,1)\ge
11n-27$.
\end{enumerate}
Combining cases 1 and 2 we obtain the estimate $C_a(T,w,1)\ge
\min[n^2, 11n-27]$. For $n\ge 8$ $n^2>11n-27$, so the inequality
$M_2(T)=C_a(T,w,1)\ge 11n-27$ holds.

For $n<8$ $n^2<11n-27$ and, thus, the optimal tree is a star
$K_{1,n}$. Let us prove that for any tree $T_{4,5}$ with $n\ge9$
pendent vertices, in which stem vertices have degrees $4$ or $5$
while the rest internal vertices having degree $3$,
$M_2(T_{4,5})=11n-27$. Consider such a tree with $s_4$ stem vertices
of degree $4$ and $s_5$ stem vertices of degree $5$. Each pendent
vertex is assigned to some stem, so the balance equation
$3s_4+4s_5=n$ holds. Note that for any $n\ge 9$ $s_4\ge 0$ and
$s_5\ge 0$ can be chosen to fulfill the balance, so the tree
$T_{4,5}$ does exist for $n\ge 9$.

The edge set $E(T_{4,5})=S_4\cup S_5\cup E_4\cup E_5\cup E_I$,
where:
\begin{itemize}
    \item $S_4$ is the set of stem edges incident to stem vertices of degree
    $4$,
    \item $S_5$ is the set of stem edges incident to stem vertices of degree $5$,
    \item $E_4$ is the set of edges connecting stem vertices of degree $4$ to internal vertices,
    \item $E_5$ is the set of edges connecting stem vertices of degree $5$ to internal
    vertices,
    \item $E_I$ is the set of edges connecting non-stem internal
    vertices.
\end{itemize}

Obviously, $|S_4|=3s_4$ and, according to (\ref{M2_def}), each edge
makes the contribution of 4 to the index $M_2(T_{4,5})$,
$|S_5|=3s_4$ and each edge from $S_5$ makes the contribution of 5.
$|E_4|=s_4$ and, as any edge from $E_4$ connects the stem vertex of
degree $4$ with an internal vertex of degree $3$, its contribution
is 12. Similarly, the contribution of each of $s_5$ edges from $E_5$
is 15.

Finally, consider a ``defoliated'' tree $T_b$ obtained from
$T_{4,5}$ by deleting all pendent vertices and stem edges. Stem
vertices of $T_{4,5}$ become leaves in $T_b$ and $E(T_b)=E_4\cup
E_5\cup E_I$. By construction, $T_b$ is a 3-tree with $s_4+s_5$
pendent vertices, so it consists of $|V(T_b)|=2s_4+2s_5-2$ vertices
and $|E(T_b)|=2s_4+2s_5-3$ edges. Thus,
$|E_I|=|E(T_b)|-|E_4|-|E_5|=s_4+s_5-3$. Each edge from $E_I$
connects two vertices of degree $3$ and, thus, makes the
contribution of 9 to $M_2(T_{4,5})$. Summing up all contributions we
have:
$$M_2(T_{4,5})=4|S_4|+5|S_5|+12|E_4|+15|S_5|+9|E_I|=33s_4+44s_5-27.$$
Taking into account the balance equation $3s_4+4s_5=n$ we finally
obtain $M_2(T_{4,5})=11n-27$ irrespective of the values of $s_4$ and
$s_5$.
\end{proof}

\begin{figure}[htpb]
\begin{center}
  \includegraphics[width=12cm]{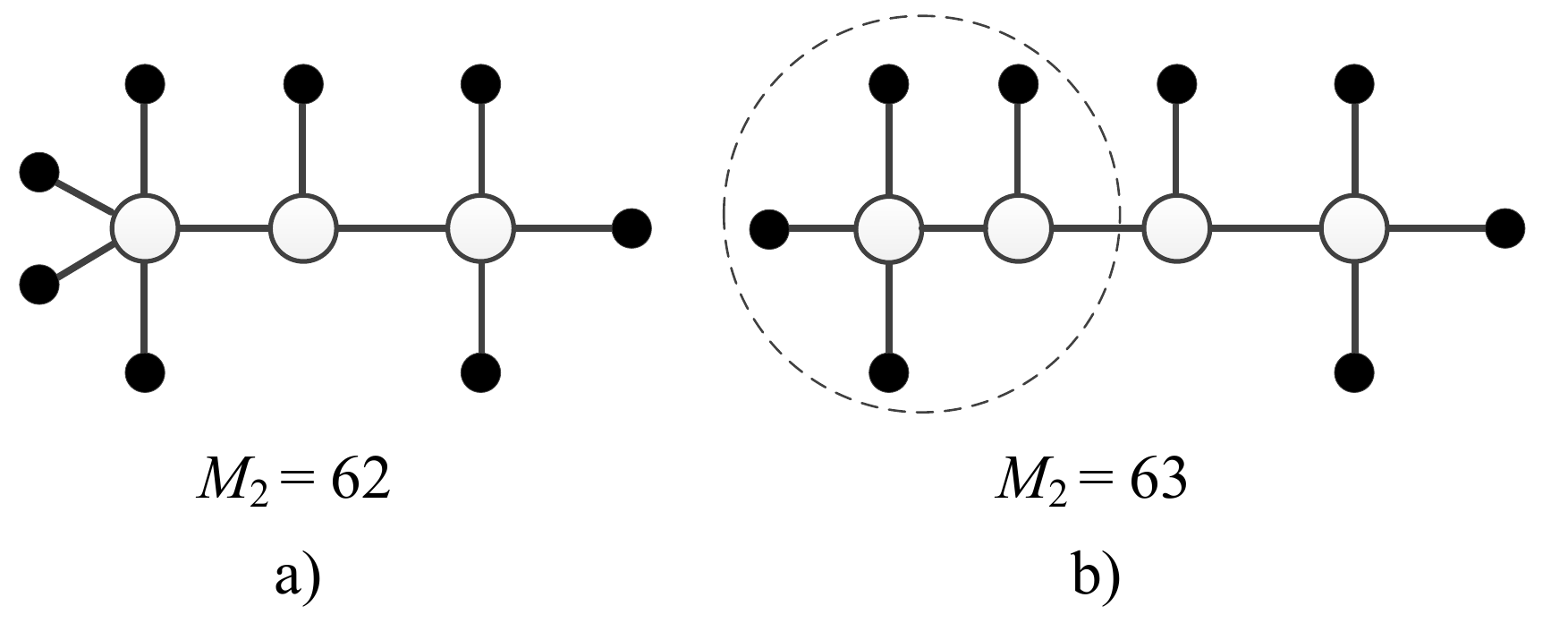}
  \caption{The $M_2$-minimal tree and the second-best tree for $n=8$}
  \label{fig_M2_8}
\end{center}
\end{figure}

Theorem \ref{M_2_theorem} provides trees which minimize $M_2$ over
all trees with $n\ge 9$ pendent vertices. From the proof of Theorem
\ref{M_2_theorem} we know that for $n<8$ the optimal tree is a star
$K_{1,n}$. The optimal tree for $n=8$ is shown in
Fig.~\ref{fig_M2_8}a.\footnote{We used a quasi-polynomial algorithm
from \cite{Goubko09} to enumerate attached trees.}

Theorem \ref{M_2_theorem} says that for some $n$ the trees, which
minimize $M_2$, are not chemical graphs. An example is shown in
Fig.~\ref{fig_M2_8}a. The optimal chemical tree for $n=8$ is shown
in Fig.~\ref{fig_M2_8}b. This tree corresponds to trans-$2$-butene
$\text{C}_4\text{H}_8$ if hydrogen atoms are not suppressed from the
diagram or to triisobutylene $\text{C}_{12}\text{H}_{24}$ otherwise.

Actually, if $n\hspace{-5pt}\mod3=1$, there should be at least one
stem vertex of degree $5$ in the optimal tree $T_{4,5}$, if
$n\hspace{-5pt}\mod3=2$, then at least two stem vertices of degree
$5$ are required to build the optimal tree $T_{4,5}$. From the proof
of Theorem \ref{M_2_theorem} one can conclude that the lower bound
$11n-27$ is not achievable with chemical graphs in these cases.

\begin{figure}[htpb]
\begin{center}
  \includegraphics[width=9cm]{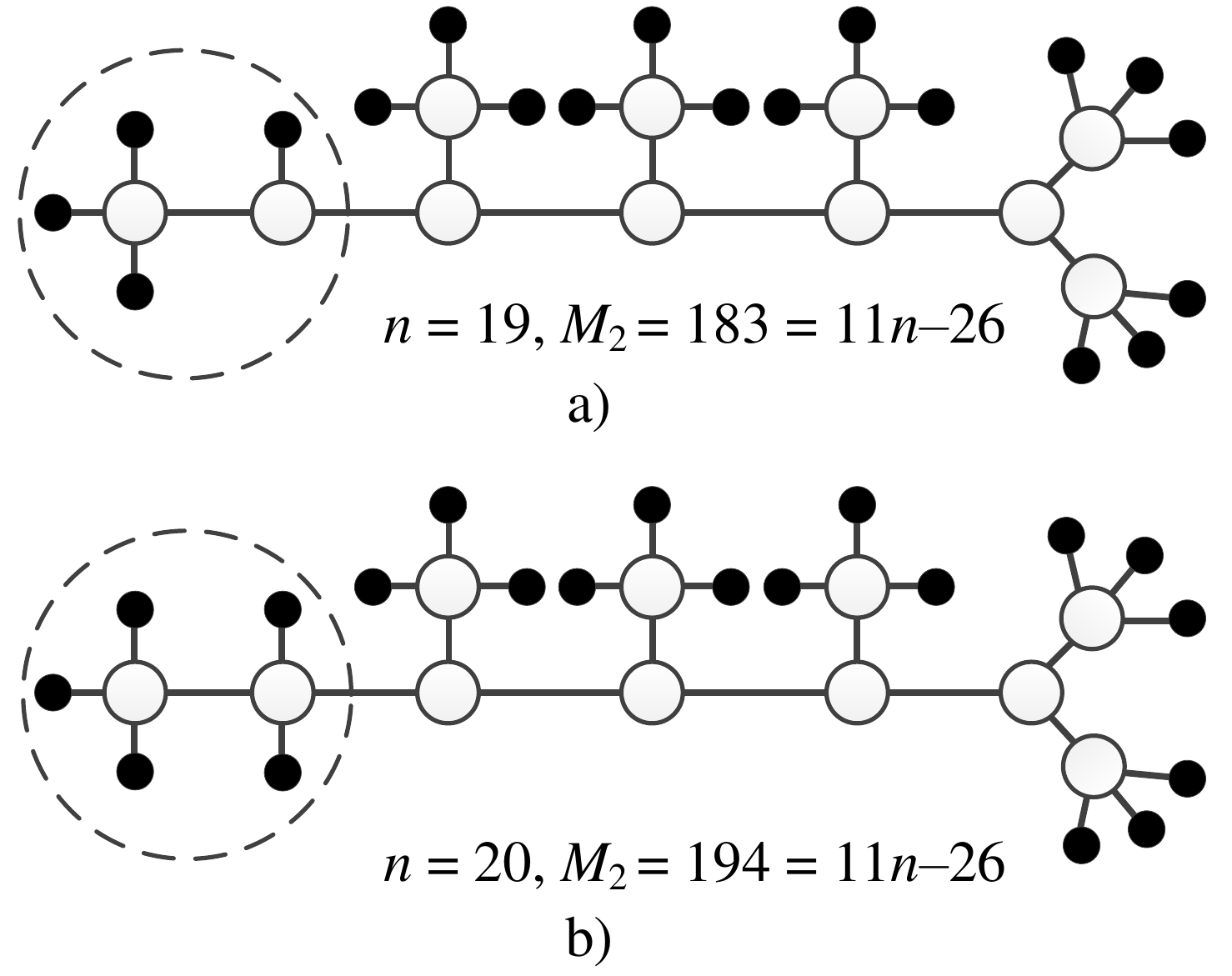}
  \caption{Examples of chemical trees minimizing $M_2$ for $n=19,20$}
  \label{fig_M2_19_20}
\end{center}
\end{figure}

At the same time, for $n\hspace{-5pt}\mod3=1$ replacement of the
subtree rooted in the stem vertex of degree $5$ with the subtree
enclosed in a dashed circle in Fig.~\ref{fig_M2_8}b gives a chemical
graph with the value $M_2$, which is only one more than the lower
bound $11n-27$. This graph appears to be the optimal chemical graph
when $n\hspace{-5pt}\mod3=1$. Analogously, for
$n\hspace{-5pt}\mod3=2$ replacement of two stem vertices and their
incident pendents with the fragment from Fig.~\ref{fig_M2_8}b gives
a chemical tree with $M_2=11n-25$, yet this tree is not the best
chemical tree for this $n$. The proof of Theorem \ref{M_2_theorem}
can be easily adopted to justify this claim (the footnote in the
proof marks the place of possible adjustment) but one can better
find a counterexample with direct enumeration of all optimal
chemical graphs with the algorithm of complexity $n^4$ from
\cite{Goubko09}. Examples of optimal chemical trees for $n=19$ and
$20$ are depicted in Fig.~\ref{fig_M2_19_20}.

\section{Conclusion}
Above we suggested an optimization framework for degree-based
indices of undirected trees. Using the discussed approach one can
calculate lower bounds for Zagreb-like indices and find the graphs
minimizing these indices over the set of trees (or chemical trees)
with the fixed number of pendent vertices.

Theorem \ref{LB_M1_theorem} provides a tight lower-bound estimate
for $M_1$ and shows that it is achieved at $4$-trees. Theorem
\ref{LB_M1_general_theorem} gives a high-quality lower-bound
estimate for the generalized $M_1$-like index. Theorem
\ref{M_2_theorem} proves the tight lower-bound estimate for $M_2$
and characterizes $M_2$-minimal trees.

Although one can surely suggest a simpler reasoning for
theorems~\ref{LB_M1_general_theorem} and \ref{M_2_theorem}, the
above proofs have an advantage, as they are open for generalization
to other degree-based graph indices, e.g., to the general Randi\'{c}
index, which is defined as $R_{\alpha}(G):=\sum_{uv\in
E(G)}d_G(u)^{\alpha}d_G(v)^{\alpha}$ (also known as $\alpha$-weight,
see \cite{Bollobas99}), or even to the abstract degree-based
topological index
$$C(G):=\sum_{v\in V(G)}c_1(d_G(v))+\sum_{uv\in
E(G)}c_2(d_G(u), d_G(v)),$$ where $c_1(d)$ is a non-negative
function of a natural argument and $c_2(d_1,d_2)$ is a non-negative
symmetric function of natural arguments. This index generalizes
almost all known topological graph indices based on vertex degrees.
As an example, one may employ the outline of the proof of Theorem
\ref{M_2_theorem} to justify the lower-bound estimate $61n/3-46$ for
the sum $M_1+M_2$. This estimate holds for trees with the number of
pendent vertices $n\ge 6$.

The proofs of theorems \ref{LB_M1_general_theorem} and
\ref{M_2_theorem}, in fact, appeal to the technique developed in
\cite{Goubko06,Goubko08,Goubko09} for directed trees with the fixed
set of leaves. As the framework developed there is not limited to
the case of degree-based topological indices, it seems promising to
apply this approach to analyze trees with the fixed number of
pendent vertices, which minimize complex topological indices:
distance-based ones (like the Wiener index), or linear combinations
of distance- and degree-based indices (some settings are provided in
\cite{LiuLiu12, Wiener47}).

\end{document}

\maketitle

\thispagestyle{empty}

\begin{abstract}
Theorem 3 in \cite{Goubko14} says that the second Zagreb index $M_2$ cannot be less than $11n-27$ for a tree with $n\ge
8$ pendent vertices. Yet the tree exists with $n=8$ vertices (the two-sided broom) violating this inequality. The
reason is that the proof of Theorem 3 relays on a tacit assumption that an index-minimizing  tree contains no vertices
of degree $2$. This assumption appears to be invalid in general. In this erratum we show that the inequality $M_2 \ge
11n-27$ still holds for trees with $n\ge 9$ vertices and provide the valid proof of the (corrected) Theorem 3.
\end{abstract}

\baselineskip=0.29in

Let $G$ be a simple connected undirected graph with the vertex set $V(G)$ and the edge set $E(G)$. Denote by $d_G(v)$
the degree of a vertex $v\in V(G)$ in the graph $G$, i.e., the number of vertices being incident to $v$ in $G$. The
second Zagreb index is defined as
\begin {equation}\label{M2_def}
M_2(G):=\sum_{uv\in E(G)} d_G(u)\,d_G(v) \ .
\end{equation}

The vertex $v\in V(G)$ with $d_G(v)=1$ is called a \emph{pendent} vertex. All others are \emph{internal} vertices. A
connected graph $T$ with $N$ vertices and $N-1$ edges is called a tree. An internal vertex in a tree is called a
\emph{stem vertex} if it has at most one incident internal vertex.

In \cite{Goubko14} the following theorem was stated.
\begin{old_theorem}\label{M_2_theorem_initial}
For any tree $T$ with $n\ge 8$ pendent vertices $M_2(T)\ge 11n-27$. The equality holds if each stem vertex in $T$ has
degree $4$ or $5$ while the other internal vertices have degree $3$. At least one such tree exits for any $n\ge 9$.
\end{old_theorem}

\begin{figure}[htpb]
\begin{center}
  \includegraphics[width=12cm]{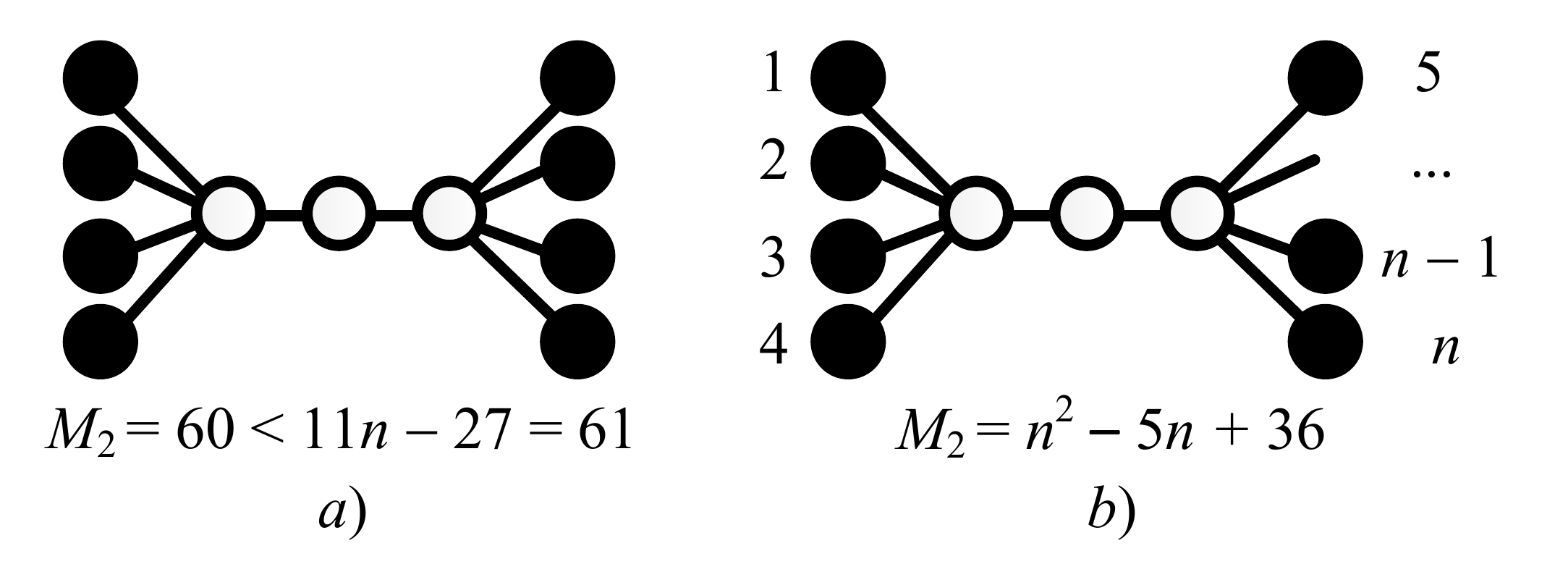}
  \caption{Two-sided brooms}
  \label{fig_broom8}
\end{center}
\end{figure}

However later a tree was found with $n=8$ pendent vertices (the two-sided broom $D(4;3;4)$, see Fig. \ref{fig_broom8}a)
with $M_2(D(4;3;4))=60<11n-27$. The inaccuracy in the proof of theorem 3 is originated from the following statement:
``...as before, we can restrict attention to the trees where all internal vertices ... have degree at least $3$...''.
This assumption is not valid in general. For example, for $8$ pendent vertices the $M_2$-minimizing tree is presented
in Fig. \ref{fig_broom8}a.

Below we show that the statement of Theorem 3 is still valid for trees with $n\ge 9$ pendent vertices by proving the
following corrected version of Theorem 3.

\begin{new_theorem}\label{M_2_theorem_corrected}
For any tree $T$ with $n\ge 9$ pendent vertices $M_2(T)\ge 11n-27$. The equality holds if each stem vertex in $T$ has
degree $4$ or $5$ while other internal vertices having degree $3$. At least one such tree exists for any $n\ge 9$.
\end{new_theorem}

We will need the following auxiliary results. Below any tree minimizing $M_2$ over the set of all trees with $n$
pendent vertices is called \emph{optimal} for short.

\begin{lemma}\label{lemma_no_incident_degree2}
For any edge $vv'\in E(T)$ in an optimal tree $T$ with $n\ge3$ pendent vertices either $d_T(v)\ge 3$ or $d_T(v')\ge 3$.
\end{lemma}
\begin{proof}
Assume that, by contradiction, $d_T(v)=d\le 2$ and $d_T(v')=d'\le 2$. Either $d=2$ or $d'=2$, as otherwise $T = K_2$
and, thus, the tree $T$ cannot have $n\ge 3$ pendent vertices. With no loss of generality suppose that $d'=2$, and,
thus, $d\le d'$. Let $v''\in V(T)$ be the second vertex incident to $v'$, and define $d'':=d_T(v'')$. Let us consider
the tree $T'$ obtained from $T$ by deleting the internal vertex $v'$ with its incident edges and adding the edge
$vv''$. Obviously, $M_2(T')-M_2(T)=dd''-dd'-d'd''$, and from $d\le d'$ we find that $M_2(T')-M_2(T)\le -dd'<0$. The
trees $T$ and $T'$ have the same number of pendent vertices, so $T$ cannot be optimal. This contradiction completes the
proof.
\end{proof}

\begin{lemma}\label{lemma_not_star}
In an optimal tree with $n\ge 8$ pendent vertices any internal vertex has at least one incident internal vertex.
\end{lemma}
\begin{proof}
If the lemma is not valid, an optimal tree is a star $K_{1,n}$ with $M_2(K_{1,n})=n^2$. Consider a two-sided broom
$D(4;3;n-4)$ (see Fig. \ref{fig_broom8}b) with $n$ pendent vertices and $M_2(D(4;3;n-4))=n^2-5n+36$. For $n\ge 8$
$M_2(D(4;3;n-4))<M_2(K_{1,n})$, so $K_{1,n}$ cannot be optimal. This contradiction completes the proof.
\end{proof}

\begin{lemma}\label{lemma_degree_le6}
Any vertex degree is at most $6$ in an optimal tree with $n\ge 8$ pendent vertices.
\end{lemma}
\begin{proof}
Assume, by contradiction, that in an optimal tree $T$ some vertex $v\in V(T)$ has degree $d_T(v)=p>6$. Let $v_1, ...,
v_p$ be its incident vertices with degrees $d_1\ge ... \ge d_p$ respectively. From Lemma \ref{lemma_not_star} we know
that $d_1\ge 2$.

Let $T'$ be a tree obtained from $T$ by adding vertices $v'$ and $v''$, edges $vv'$ and $v'v''$, and redirecting edges
$vv_i$, $i=4,...,p$, to the vertex $v''$ instead of $v$ (see Fig. \ref{fig_degree6}).
\begin{figure}[htpb]
\begin{center}
  \includegraphics[width=12cm]{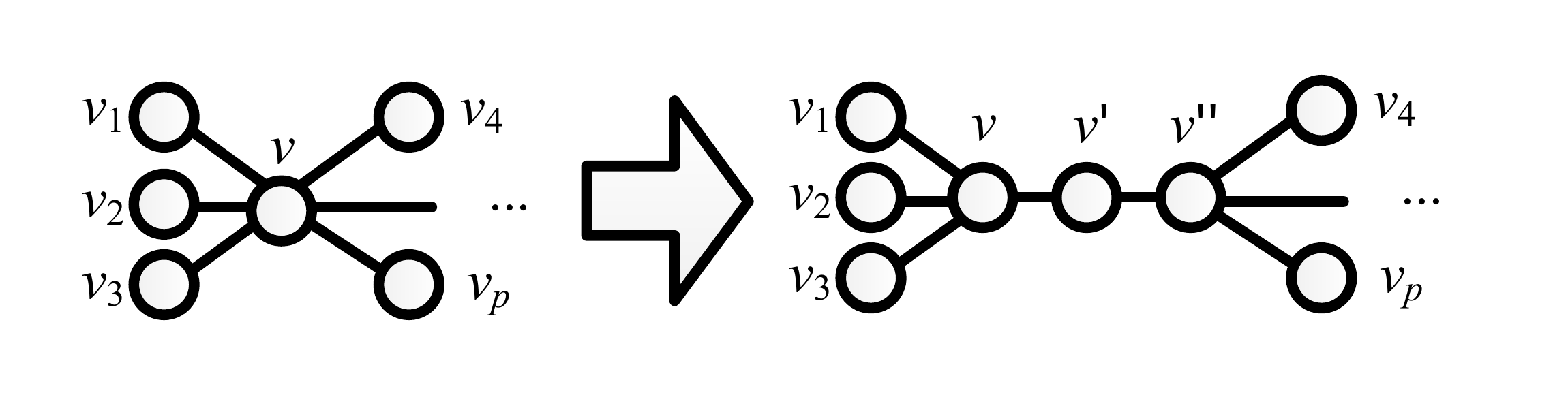}
  \caption{Transformation of the vertex with degree $d\ge 6$}
  \label{fig_degree6}
\end{center}
\end{figure}
\begin{multline}\label{delta_6}
\Delta:=M_2(T')-M_2(T)=\sum_{i=1}^3 4\,d_i+2\cdot4+2\,(p-2)+\sum_{i=4}^p (p-2)\,d_i-\sum_{i=1}^p pd_i=\\
=2\,p+4-2\sum_{i=4}^p d_i-(p-4)\sum_{i=1}^3 d_i.
\end{multline}
If $p\ge7$ then $p-4>0$ in (\ref{delta_6}). From $d_1\ge 2$, $d_i\ge 1$, $i=2,...,p$, it follows that $\sum_{i=1}^3
d_i\ge 4$. Thus, $\Delta\le 2\,p+4-2\,(p-3)-4\,(p-4)=26-4\,p<0$. The trees $T$ and $T'$ have the same number of pendent
vertices, so $T$ cannot be optimal. This contradiction completes the proof.
\end{proof}

An \emph{attached tree} is a rooted tree with a root being a pendent vertex parameterized with some ``virtual degree''
(degree of the vertex this tree is ``attached'' to). The vertex incident to the root is called a \emph{sub-root}, and
if the sub-root has degree 2, its incident vertex other than root is called a \emph{sub-sub-root}. It will be
convenient to consider the root as a non-pendent vertex.

The \emph{cost} of an attached tree $T$ with some root $w$ of ``virtual degree'' $p$ and a sub-root $m$ of degree $d$
is defined as
\begin{equation}\label{M2_attached_tree_cost}
C_a(T,w,p):=p\,d+\sum_{uv\in E(T)\backslash \{wm\}}d_T(u)\,d_T(v).
\end{equation}

Consider a tree $T$ and fix any vertex $v\in V(T)$. If it has degree $p$ and incident vertices $v_1, ..., v_p$, then
$T$ is a union of $p$ attached trees $T_1,..., T_p$ with the common root $v$ and sub-roots $v_1, ..., v_p$, so
$M_2(T)=\sum_{i=1}^p C_a(T_i,v, p)$. Below we limit attention to the attached trees, which can be met in optimal trees,
so Lemmas \ref{lemma_no_incident_degree2}, \ref{lemma_not_star}, and \ref{lemma_degree_le6} are supposed to be valid
for every attached tree in hand.

Let $\mathcal{T}_a(n,p)$ be the collection of all attached trees with $n$ pendent vertices where the root (denoted with
$w$) has degree $p\ge 2$ (remember the root is considered non-pendent), and introduce the cost of an optimal attached
tree $C_a^*(n,p)=\min_{T\in \mathcal{T}_a(n,p)}C_a(T,w,p)$.

\begin{lemma}\label{lemma_Bellman}
For $n=1$ $C_a^*(n,p)=p$, while for $n\ge 2$
\begin{equation}\label{bellman_equation}
C_a^*(n,p)= \min_{d=2,...,6}\min_{n_1,...,n_{d-1}}\left\{pd+\sum_{i=1}^{d-1} C_a^*(n_i,d)|n_i\ge1,
\sum_{i=1}^{d-1}n_i=n\right\}.
\end{equation}
\end{lemma}
\begin{proof}
The case of $n=1$ is obvious. For $n\ge2$ each combination of $d$ and $n_1,...,n_{d-1}$ in the right-hand side of
(\ref{bellman_equation}) gives rise to an attached tree with $n$ pendent vertices, the sub-root enjoying degree $d$ and
being a root of $d-1$ optimal attached trees with $n_1, ..., n_{d-1}$ pendent vertices respectively. So, $C_a^*(n,p)$
cannot exceed the right-hand side in (\ref{bellman_equation}).

Let an optimal attached tree $T$ with $n$ pendent vertices have a root $w$ and a sub-root $m$ of some degree $d^*$.
Then $T$ is a union of the edge $wm$ and $d^*-1\ge1$ attached sub-trees $T_1,..., T_{d^*-1}$ with the common root $m$.
Let the trees $T_1,..., T_{d^*-1}$ have $n_1^*,...,n_{d^*-1}^*$ pendent vertices respectively. By definition,
$C_a(T_i,m,d^*)\ge C_a(n_i^*,d^*)$. So,
$$C_a^*(n,p)= pd^*+\sum_{i=1}^{d^*-1} C_a(T_i^*,m, d^*)\ge pd^*+\sum_{i=1}^{d^*-1} C_a^*(n_i^*,d^*),$$
which is obviously not less than the right-hand side in (\ref{bellman_equation}).
\end{proof}

Let us rewrite (\ref{bellman_equation}) as $C_a^*(n,p)=\min\left[C_{>2}(n,p),C_2(n,p)\right]$, where

\begin{equation}\label{bellman_equation3}
C_{>2}(n,p):=\min_{d=3,...,6}\min_{n_1,...,n_{d-1}}\left\{pd+\sum_{i=1}^{d-1} C_a^*(n_i,d)|n_i\ge1,
\sum_{i=1}^{d-1}n_i=n\right\},
\end{equation}
\begin{equation}\label{bellman_equation2}
C_{2}(n,p):=2p+C_a^*(n,2).
\end{equation}

From Lemma \ref{lemma_no_incident_degree2}, vertices of degree 2 cannot be incident in an optimal tree. So, if the root
has degree 2 in an optimal attached tree, the sub-root must have some degree $d\ge3$, and, thus,
$C_a^*(n,2)=C_{>2}(n,2)$.

Now we are ready to prove Theorem $3^*$.
\begin{proof}
Let us prove the following estimate for the cost of an optimal attached tree:
\begin{equation}\label{M2_attached_inequality_new}
C_a^*(n,p)\ge \underline{C}_a(n,p)=\begin{cases} p & \mbox{if } n=1, \\
40 & \mbox{if } p=5\mbox{ and }n=4,\\
32 & \mbox{if } p=6\mbox{ and }n=3,\\
42 & \mbox{if } p=6\mbox{ and }n=4,\\
54 & \mbox{if } p=6\mbox{ and }n=5,\\
11n+3p-18 & \mbox{otherwise}, \end{cases}
\end{equation}
which is valid for $p=3,...,6$.

For $n=1$ inequality (\ref{M2_attached_inequality_new}) trivially holds. Assume it holds for all $n'< n$. Let us prove
that it also holds for $n$. As $d\ge3$ in (\ref{bellman_equation3}), so $n_i<n$, $i=1,...,d-1$, and, by induction
hypothesis, $C_a^*(n_i,d)\ge \underline{C}_a(n_i,d)$ in (\ref{bellman_equation3}). Thus, from
(\ref{bellman_equation3}), (\ref{bellman_equation2}), we are enough to prove that
\begin{equation}\label{ineq_case3}
\min_{d=3,...,6}\min_{n_1,...,n_{d-1}}\left\{pd+\sum_{i=1}^{d-1} \underline{C}_a(n_i,d)|n_i\ge1,
\sum_{i=1}^{d-1}n_i=n\right\}\ge\underline{C}_a(n,p),
\end{equation}
and
\begin{equation}\label{ineq_case2}
2p+\min_{d=3,...,6}\min_{n_1,...,n_{d-1}}\left\{2d+\sum_{i=1}^{d-1} \underline{C}_a(n_i,d)|n_i\ge1,
\sum_{i=1}^{d-1}n_i=n\right\}\ge\underline{C}_a(n,p).
\end{equation}

For $n = 2, ..., 25$ inequalities (\ref{ineq_case3}) and (\ref{ineq_case2}) are validated by exhaustive enumeration in
their left sides. Lines 2-5 in (\ref{M2_attached_inequality_new}) correspond to the situations when the optimal
attached tree is a \emph{broom} $B(3,n)$, in which a sub-root has degree 2.

Consider $n \ge 26$ and an arbitrary vector $(n_1, ..., n_{d-1})$ such that $n_i\ge 1, \sum_{i=1}^{d-1}n_i=n$. Let us
define $\delta_i:=\#\left\{j: n_j=i, j=1,...,d-1\right\}$ and $I_d:=\{i:\underline{C}_a(i,d)\neq 11n+3d-18\}$. For
$n\ge26$ at least one $\delta_i\ge1$ for $i\ge6$, so $\sum_{i\in I_d}\delta_i\le d-2$, and (\ref{ineq_case3}) can be
written as
\begin{multline}\label{ineq_case3bis}
\min_{d=3,...,6}\min_{\delta_i:i\in I_d}\{pd+\sum_{i\in I_d} \delta_i\underline{C}_a(i,d)+11\left(n-\sum_{i\in
I_d}i\delta_i\right) +\\
+(3d-11)\left(d-1-\sum_{i\in I_d}i \right)|\delta_i\ge0, \sum_{i\in I_d}\delta_i\le d-2\}\ge\underline{C}_a(n,p).
\end{multline}

For (\ref{ineq_case2}) we can perform a similar transformation.

\begin{enumerate}
\item Consider $d=3$. As $I_3=\{1\}$, substituting (\ref{M2_attached_inequality_new}) into the left-hand side of (\ref{ineq_case3bis}) obtain
\begin{gather}\min_{\delta_1=0,1}\{pd+\delta_1 d +11(n-\delta_1)+(d-1-\delta_1)(3d-18)\}=\nonumber\\
\min_{\delta_1=0,1}\{11n+3p+\delta_1-18\}=11n+3p-18,\nonumber\end{gather}
which is equal to $\underline{C}_a(n,p)$, so
(\ref{ineq_case3}) holds.

In the same way, substituting (\ref{M2_attached_inequality_new}) into the left-hand side of (\ref{ineq_case2}) we
obtain
$$\min_{\delta_1=0,1}\{2p+2d+\delta_1 d+11(n-\delta_1)+(d-1-\delta_1)(3d-18)\}=11n+2p-12,$$
which is not less than $\underline{C}_a(n,p)=11n+3p-18$ for $p\le6$, so (\ref{ineq_case2}) also holds.

\item For $d=4$ $I_4=\{1\}$, and substitution of (\ref{M2_attached_inequality_new}) into the left-hand side of (\ref{ineq_case3bis}) gives
\begin{gather}
\min_{\delta_1=0,1,2}\{pd+\delta_1 d +11(n-\delta_1)+(d-1-\delta_1)(3d-18)\}=\nonumber\\
=\min_{\delta_1=0,1,2}\{11n+4p-\delta_1-18\}=11n+4p-20,\nonumber
\end{gather}
which is greater than
$\underline{C}_a(n,p)=11n+3p-18$ for $p\ge3$, so (\ref{ineq_case3}) holds.

In the same way for (\ref{ineq_case2}) we obtain
$$\min_{\delta_1=0,1,2}\{2p+2d+\delta_1 d+11(n-\delta_1)+(d-1-\delta_1)(3d-18)\}=11n+2p-12,$$
which is not less than $\underline{C}_a(n,p)=11n+3p-18$ for $p\le6$, so (\ref{ineq_case2}) holds.

\item For $d=5$ $I_5=\{1,4\}$, so, substituting (\ref{M2_attached_inequality_new}) into the left-hand side of (\ref{ineq_case3bis}),
we obtain
\begin{gather} \min_{\delta_1, \delta_4\ge0}\{pd+\delta_1 d +40\delta_4 +
11(n-\delta_1-4\delta_4)+(d-1-\delta_1-\delta_4)(3d-18)|\nonumber\\
|\delta_1+ \delta_4 \le d-2\}=\min_{\delta_1, \delta_4\ge0}\{11n+5p-3\delta_1-\delta_4-12|\delta_1+ \delta_4 \le
3\}.\nonumber
\end{gather}
The
latter minimum is attained at $\delta_1=3,\delta_4=0$ and is equal to $11n+5p-21$, which is greater than
$\underline{C}_a(n,p)=11n+3p-18$ for $p\ge2$, so (\ref{ineq_case3}) holds.

Analogously for (\ref{ineq_case2}) we obtain \begin{gather}\min_{\delta_1, \delta_4\ge0}\{2p+2d+\delta_1 d +40\delta_4
+ 11(n-\delta_1-4\delta_4)+(d-1-\delta_1-\delta_4)(3d-18)|\nonumber\\
|\delta_1+ \delta_4 \le d-2\}=\min_{\delta_1, \delta_4\ge0}\{11n+2p-3\delta_1-\delta_4-2|\delta_1+ \delta_4 \le
3\}=11n+2p-11,\nonumber\end{gather} which is always greater than $\underline{C}_a(n,p)=11n+3p-18$ for $p\le6$, so
(\ref{ineq_case2}) holds.

\item For $d=6$ $I_6=\{1,3,4,5\}$, so we write the left-hand side of (\ref{ineq_case3bis}) as
\begin{gather}
\min_{\delta_1, \delta_3, \delta_4,\delta_5\ge0}\{pd+\delta_1 d +32\delta_3+42\delta_4+54\delta_5 +
11(n-\delta_1-3\delta_3-4\delta_4-5\delta_5)+\nonumber\\
+(d-1-\delta_1-\delta_3-\delta_4-\delta_5)(3d-18)|\delta_1+
\delta_3+\delta_4+\delta_5 \le d-2\}=\nonumber\\
=\min_{\delta_1, \delta_3, \delta_4,\delta_5\ge0}\{11n+6p-5\delta_1-\delta_3-2\delta_4-\delta_5|\delta_1+
\delta_3+\delta_4+\delta_5 \le 4\}.\nonumber
\end{gather}

The minimum is attained at $\delta_1=4,\delta_3=\delta_4=\delta_5=0$ and is equal to $11n+6p-20$, which is greater than
$\underline{C}_a(n,p)=11n+3p-18$ for $p\ge 1$, so (\ref{ineq_case3}) holds.

In the same way we prove that the left-hand side in (\ref{ineq_case2}) is equal to $11n+2p-8$, which is greater than
$\underline{C}_a(n,p)=11n+3p-18$ for $p\le6$, so (\ref{ineq_case2}) holds.
%In the same way we transform the left-hand side of (\ref{ineq_case2}) into
%\begin{gather}
%\min_{\delta_1, \delta_3, \delta_4,\delta_5\ge0}\{2p+ 2d +\delta_1 d +32\delta_3+42\delta_4+54\delta_5 +
%11(n-\delta_1-3\delta_3-4\delta_4-5\delta_5)+\nonumber\\
%+(d-1-\delta_1-\delta_3-\delta_4-\delta_5)(3d-18)|\delta_1+ \delta_3+\delta_4+\delta_5 \le d-2\}=\nonumber\\
%\min_{\delta_1, \delta_3, \delta_4,\delta_5\ge0}\{11n+2p+12-5\delta_1-\delta_3-2\delta_4-\delta_5|\delta_1+
%\delta_3+\delta_4+\delta_5 \le 4\}.\nonumber
%\end{gather}
%The minimum is attained at $\delta_1=4,\delta_3=\delta_4=\delta_5=0$ and is equal to $11n+2p-8$, which is greater than
%$\underline{C}_a(n,p)=11n+3p-18$ for $p\le6$, so (\ref{ineq_case2}) holds.
\end{enumerate}
Thus, we proved inequality (\ref{M2_attached_inequality_new}). Let us prove that $M_2(T)\ge 11n-27$ for every tree with
$n$ pendent vertices. By Lemma \ref{lemma_not_star}, the optimal tree $T$ is not a star for $n\ge9$, so every internal
vertex in an optimal tree has an incident internal vertex. At least one of them must be a stem vertex $m$ of some
degree $d$, which has has $d-1$ incident pendent vertices and one incident internal vertex. From Lemmas
\ref{lemma_no_incident_degree2} and \ref{lemma_degree_le6}, $d\in \{3,...,6\}$.

So, the value of the index $M_2(T)$ adds up from the total contribution $(d-1)d$ of $d-1$ pendent vertices and the cost
of an attached sub-tree $T_1$:
$$M_2(T)=(d-1)d+C_a(T_1,m,d)\ge (d-1)d+\underline{C}_a(n-d+1,d).$$

Consider $n\ge11$, so that $n-d+1\ge6$ and $\underline{C}_a(n-d+1,d)$ is always equal to $11(n-d+1)+3d-18$. In this
case $M_2(T)\ge 11n+(d-9)d-7$. The minimum in the right-hand side is attained at $d=4,5$ and is equal to $11n-27$.

For $n=9,10$ we also need to check that double brooms $D(d-1,3,n-d+1)$, originated from the lines 2-5 in
(\ref{M2_attached_inequality_new}), do not violate inequality $M_2(D(d-1,3,n-d+1))\ge 11n-27$. For example,
$M_2(D(5,3,4))\ge 5\cdot6+C_a(4,6)=11n-27=72$, $M_2(D(5,3,5))\ge 5\cdot6+C_a(5,6)=84>11n-27=83$.

The existence of the optimal tree is proved as in Theorem 3 in \cite{Goubko14}.
\end{proof}